\documentclass[12pt,reqno]{amsart}

\usepackage{latexsym}
\usepackage{amsmath}
\usepackage{amscd}
\usepackage{amssymb}
\usepackage{mathrsfs}
\usepackage{comment}
\usepackage{color}
\usepackage{enumerate}
\usepackage{soul}

\usepackage{graphicx}
\usepackage{subfigure}
\usepackage{url}

\newtheorem{theorem}{Theorem}[section]

\numberwithin{equation}{section}

\begin{document}
\title{Noise and dissipation in rigid body motion}
\author[A. Arnaudon]{Alexis Arnaudon}
\author[A. L. De Castro]{Alex L. De Castro}
\author[D. D. Holm]{Darryl D. Holm}

\address{AA, DH: Department of Mathematics, Imperial College, London SW7 2AZ, UK}
\address{AC: Department of Mathematics, Imperial College, London SW7 2AZ, UK and Departamento de Matem\'atica PUC-Rio, Rio de Janeiro 22451-900.}

\maketitle

\begin{abstract}
Using the rigid body as an example, we illustrate some features of stochastic geometric mechanics. These features include: i) a geometric variational motivation for the noise structure involving Lie-Poisson brackets and momentum maps, ii) stochastic coadjoint motion with double bracket dissipation, iii) the Lie-Poisson Fokker-Planck description and its stationary solutions, iv) random dynamical systems, random attractors and SRB measures  connected to statistical physics.
\end{abstract}

\setcounter{tocdepth}{1}

\section{Introduction}

The rigid body sets the paradigm for geometric mechanics. Any new ideas in this field must always be tested on the rigid body. Therefore, to illustrate the effects of stochasticity in geometric mechanics we may begin with the rigid body. The key idea underlying geometric mechanics is coadjoint motion on level sets of momentum maps, derived via reduction by symmetry as an equivariant map from the original phase space to the dual of a Lie algebra of the symmetry group of the Hamiltonian, or Lagrangian. The resulting Hamiltonian formulation involves a Lie-Poisson bracket defined on the dual Lie algebra of the symmetry. Since the Lie-Poisson bracket preserves coadjoint orbits, one may incorporate dissipation as a double Lie-Poisson bracket, which is compatible with coadjoint motion. The probability distribution on the coadjoint orbit for the dynamics of the stochastic rigid body without double bracket dissipation is described by the Lie-Poisson Fokker-Planck equation, whose asymptotic solution tends to a constant on each coadjoint orbit. However, when the double bracket dissipation is included, the nonlinear interaction between noise and dissipation leads to an equilibrium probability distribution given by a Gibbs measure of the standard form $\mathbb{P}_\infty(\mathbf \Pi)= Z^{-1}\exp(-E(\mathbf \Pi)/k_BT)$, which depends on the energy $E(\mathbf \Pi)$ as a function of body angular momentum $\mathbf \Pi$,  as well as constants $k_B$ (Boltzmann constant), $T$ (temperature), and the normalising factor $Z$. The time-dependent approach to equilibrium with double bracket dissipation creates SRB (Sinai-Ruelle-Bowen) measures. Thus, introducing  double bracket dissipation allows us to study random attractors on level sets of coadjoint orbits for stochastic rigid body dynamics. For more details about stochastic coadjoint motion, existence of SRB measures and extension to semidirect product theory with an application to the heavy top, we refer to \cite{arnaudon2016noise}.

\section{Theory}\label{noise}

\subsection{Preliminaries}

This paper focuses on dynamical systems of rigid body type, which are written on the dual of a semi-simple Lie algebra, such as $\mathfrak{so}(3)$ in the case of the classical rigid body.  
The property of semi-simplicity guarantees the existence of a non-degenerate pairing, called a Killing form, which allows the identification of the Lie algebra with its dual. 
We will denote this pairing by $\langle \xi, \eta\rangle := \mathrm{Tr}(\mathrm{ad}_\xi\mathrm{ad}_\eta)= \epsilon \mathrm{Tr}(\xi\eta)$, for $\xi,\eta\in \mathfrak g$ and a number $\epsilon$ which depends on the Lie algebra. 
This pairing is bi-invariant, i.e., $\kappa (\xi, \mathrm{ad}_\zeta \eta) = \kappa (\mathrm{ad}_\xi \zeta,\eta)$ and furthermore, if the Lie algebra is compact, $\epsilon<0$ and (minus) the Killing form defines a norm. 
The Killing form will allow us to reformulate coadjoint equations on the dual Lie algebra, in terms of ad-operations (Lie brackets) on the Lie algebra. (See \cite{arnaudon2016noise} for the corresponding equations written on the dual Lie algebra.)
We will find that all of our equations of motion will contain a Lie bracket, which of course depends on the Lie algebra. These equations characterise coadjoint motion which is always restricted to some particular submanifolds called coadjoint orbits, parameterised by the initial conditions. In the free rigid body example, the coadjoint orbits are spheres in the three-dimensional space of angular momenta. 

Here, we will only use the simplest type of noise, composed of $n$ independent Wiener processes $W_t^i$ indexed by $i=1,2,\dots,n$. See for example \cite{chen2015constrained, ikeda2014stochastic} for more details about stochastic processes. 
For simplicity, $n$ will be the dimension of the Lie algebra, which comprises the dynamical variables of our system; although, in principle, $n$ could be arbitrary.
In the stochastic integrals we discuss here, we will use the multiplication symbol $(\,\circ\,)$, which denotes a stochastic integral in the Stratonovich sense. While the Stratonovich integral admits the normal rules of calculus, the It\^o integral requires the It\^o calculus in the computations. Both representations are equivalent, provided a correction term is added, but the Stratonovich integral will be more convenient for us in dealing with variational calculus. Our strategy is to stay with the Stratonovich sense until nearly the end of our calculations, then transform to the It\^o representation when passing to the Fokker-Planck equation. 

\subsection{Structure preserving stochastic deformations}\label{EP-noise-section}
A natural framework for understanding dynamical systems with symmetry is Lie group reduction \cite{marsden1999intro}. This approach leads to reduced equations in terms of an equivariant momentum map taking values on the dual Lie algebra of the Lie symmetry group, and evolving by coadjoint motion. However, here we will use a different approach, which yields equivalent dynamical equations for the same momentum map, while avoiding the use of Lie groups. 
Thus, stochasticity on Lie groups will not be discussed here. (See \cite{arnaudon2014stochastic}, for discussions of Lie group reduction by symmetry for stochastic variational principles.)
The equivalent formulation we discuss is the so-called Clebsch principle, which constrains the variations to respect certain auxiliary evolution equations for a set of configuration variables $q$ constrained by Lagrange multipliers $p$. These are the Clebsch variables $(q,p)$, which in the present case we take to be $(q,p)\in \mathfrak g\times \mathfrak g$, after using the Killing form to identify $\mathfrak g^*$ with $\mathfrak g$ for the $p$ variables. 
This method is often used in the design of control systems, where our dynamical variable $\xi\in \mathfrak g$ below would play the role of a control parameter. 

To be concise, we combine the Clebsch variational principle with noise. 
For this, we first introduce the so-called stochastic potentials $\Phi_i: \mathfrak g\to \mathbb R$ which are prescribed functions of the momentum map $\mu:= \frac{\partial l(\xi)}{\partial \xi}={\rm ad}_qp$, for the reduced Lagrangian $l(\xi)$ whose Hamilton's principle governs the deterministic system for $\xi$. 
Then, we write the following constrained stochastic variational principle,
\begin{align}
	S(\xi,q,p) =\int  l(\xi) dt + \int \langle p, dq + \mathrm{ad}_\xi q\, dt \rangle   + \int \sum_{i=1}^n \Phi_i(\mu )\circ dW^i_t.
		\label{Sto-Clebsch-action}
\end{align}
Following the detailed calculations in \cite{holm2015variational}, the free variations of this action functional yield the stochastic Euler-Poincar\'e equation,
\begin{align}
	d \frac{\partial l(\xi)}{\partial \xi} + \mathrm{ad}_\xi \frac{\partial l(\xi)}{\partial \xi} dt 
	- \sum_i \mathrm{ad}_\frac{\partial \Phi_i(\mu)}{\partial \mu}  \frac{\partial l(\xi)}{\partial \xi}\circ dW_t^i=0\,.
	\label{SEP}
\end{align}
After having defined the Stratonovich stochastic process \eqref{SEP}, one may compute its corresponding It\^o form.
For convenience, we denote the field  $\sigma_i:= -\,\frac{\partial \Phi_i(\mu)}{\partial \mu} \in \mathfrak g$. A direct calculation gives the It\^o correction  
$- \frac12 \sum_i \mathrm{ad}_{\sigma_i}\left ( \mathrm{ad}_{\sigma_i}\frac{\partial l(\xi)}{\partial \xi}\right ) dt$ which must be added to \eqref{SEP} in order to interpret the stochastic integral as an It\^o integral. 
Let us simplify matters further by rewriting the It\^o version of \eqref{SEP} in term of $\mu:= \frac{\partial l(\xi)}{\partial \xi} $ as
\begin{align}
	d \mu + \mathrm{ad}_\xi \mu \mbox{ }dt + \sum_i \mathrm{ad}_{\sigma_i}  \mu  \mbox{ }dW_t^i 
	- \frac12 \sum_i \mathrm{ad}_{\sigma_i}\left (\mathrm{ad}_{\sigma_i}\mu\right ) dt=0\,.
	\label{SEP-simple-ito-mu}
\end{align}
The Lie-Poisson formulation of this equation results: 
\begin{align}
\begin{split}
	d f(\mu) &= \left \langle \mu, \left[ \frac{\partial f}{\partial \mu},  \frac{\partial h}{\partial \mu}\right ] \right \rangle dt + \sum_i\left \langle \mu, \left[ \frac{\partial f}{\partial \mu},  \frac{\partial \Phi_i}{\partial \mu}\right ] \right \rangle \circ dW_t^i\\
	&=: \{f,h\} dt +\sum_i \{ f,\Phi_i\} \circ dW_i\,,
\\
	&= \{f,h\} dt +\sum_i \{ f,\Phi_i\}  dW_i 
	+ \frac12 \sum_i\{ \{ f,\Phi_i\}, \Phi_i \}\,dt
	\,,
\end{split}
\label{LP-sto}
\end{align}
where the Lie-Poisson bracket $\{\cdot ,\cdot \}$ is defined as in the deterministic case, so the motion takes place on coadjoint orbits, characterised as level sets of \emph{Casimir functions}, which comprise the kernel of the Lie-Poisson bracket defined in \eqref{LP-sto}. We refer to  \cite{arnaudon2016noise,cruzeiro2016stochastic,gaybalmaz2016geometric} for more details. 

\subsection{The Fokker-Planck equation and invariant distributions}\label{FP-section}

Next, we derive a geometric version of the classical Fokker-Planck equation using our SDE \eqref{SEP}. 
Recall that the Fokker-Planck equation describes the time evolution of the probability distribution $\mathbb{P}$ for the process driven by \eqref{SEP}. See, for example, \cite{ikeda2014stochastic} for the standard reference for stochastic processes. 
We will consider $\mathbb P$ as a normalised function on $\mathfrak g$ with values in $\mathbb R$ .
First, the generator of the process \eqref{SEP} can be readily found from its Lie-Poisson form \eqref{LP-sto} 
\begin{align}
	Lf(\mu) = \left \langle \mathrm{ad}_\xi\mu ,\frac{\partial f}{\partial \mu}\right \rangle - \sum_i\left \langle \mathrm{ad}_{\sigma_i}\mu,\frac{\partial}{\partial \mu} \left \langle \mathrm{ad}_{\sigma_i}\mu ,\frac{\partial f}{\partial \mu}\right \rangle\right \rangle,
	\label{FP-gen}
\end{align}
where $f:\mathfrak g\to \mathbb R$ is an arbitrary function of $\mu$. 
Then, provided that the  $\Phi_i$'s are linear functions of the momentum $\mu$, the diffusion terms of the infinitesimal generator $L$ will be \emph{self-adjoint} with respect to the $L^2$ pairing 
$\langle\,f\,,\, \mathbb P \,\rangle_{L^2} := \int_\mathfrak{g} f(\mu)\mathbb P(\mu)d\mu$.

The Fokker-Planck equation, $\frac{d\mathbb P}{dt}=-L\mathbb P(\mu)$, describes the dynamics of the probability distribution $\mathbb P$ associated to the stochastic process for $\mu$, in the standard advection diffusion form. 
The underlying geometry of the Fokker-Planck equation may be highlighted by rewriting it in terms of the Lie-Poisson bracket structure, as 
\begin{align}
	\frac{d}{dt} \mathbb P =-L\mathbb P(\mu) 
	= - \{ h,\mathbb P\} +  \sum_i\{\Phi_i, \{\Phi_i,\mathbb P\}\} \,,
	\label{HamFP}
\end{align}
where $h(\mu)$ is the Hamiltonian associated to $l(\xi)$ by the Legendre transform. 
In \eqref{HamFP}, we have recovered the Lie-Poisson formulation \eqref{LP-sto} of the Euler-Poincar\'e equation together with a dissipative term arising from the noise of the original SDE in a double Lie-Poisson bracket form; see  \cite{arnaudon2016noise,cruzeiro2016stochastic,gaybalmaz2016geometric} for more details. 

This formulation gives the following theorem for invariant distributions of \eqref{FP-gen}.
\begin{theorem}\label{limit-thm}
	The invariant distribution $\mathbb P_\infty$ of the Fokker-Planck equation \eqref{FP-gen}, i.e, $L\mathbb P_\infty=0$ is \emph{uniform} on the coadjoint orbits on which the SDE \eqref{SEP} evolves. 
\end{theorem}
The proof of this theorem is based on the hypo-coercive property of the Fokker-Planck operator in \eqref{HamFP}; see \cite{arnaudon2016noise}.  
Compactness of the coadjoint orbits (Casimir level sets) is necessary for the existence of a non-vanishing invariant measure. 

\subsection{Double bracket dissipation}\label{dissipation-section}

We can now add dissipation in our systems, for which the solutions of the stochastic process will still lie on the deterministic coadjoint orbit.
For this purpose, we will use \emph{double bracket dissipation}, which was studied in detail in \cite{bloch1996euler} and was generalised recently in \cite{gaybalmaz2013selective, gaybalmaz2014geometric}.

For the stochastic process \eqref{SEP}, the dissipative stochastic Euler-Poincar\'e equation written in Hamiltonian form with double bracket dissipation is 
\begin{align}		
	d\mu &+ \mathrm{ad}_\frac{\partial h}{\partial \mu} \mu\, dt 
	+ \theta\, \left [ \frac{\partial C}{\partial \mu}, \left [ \frac{\partial C}{\partial \mu}, \frac{\partial h}{\partial \mu} \right ]\right ] dt + \sum_i\mathrm{ad}_{\sigma_i} \mu \circ dW_t^i = 0 \,,
	\label{SEP-Diss}
\end{align}
where $\theta>0$ parametrises the rate of energy dissipation and $C$ is a chosen Casimir of the coadjoint orbit, i.e. an constant function on the space of solutions of the original equation.  

As before, we compute the Fokker-Planck equation for the Euler-Poincar\'e stochastic process \eqref{SEP-Diss} which is now modified by the double bracket dissipative term
\begin{align}
	\frac{d}{dt} \mathbb P(\mu) + \{h,\mathbb P\}  +\theta\left \langle\left [\frac{\partial \mathbb P}{\partial \mu}, \frac{\partial C}{\partial \mu}\right], \left [ \frac{\partial h}{\partial \mu}, \frac{\partial C}{\partial \mu}\right]\right \rangle - \frac12 \sum_i  \{\Phi_i,\{\Phi_i,\mathbb P\}\}=0.
	\label{FP-Diss}
\end{align}
The invariant distribution of this Fokker-Planck equation is no longer a constant on the coadjoint orbits. Instead, it now depends on the energy, as summarized in the following theorem. 
\begin{theorem}\label{FP-diss-thm}
	Let the noise amplitude be of the form $\sigma_i= \sigma e_i$ for an arbitrary $\sigma \in \mathbb R$, where the $e_i$'s span the underlying vector space of the Lie algebra $\mathfrak g$. 
	The invariant distribution of the Fokker-Planck equation \eqref{FP-Diss} associated to \eqref{SEP-Diss} with Casimir $C= \kappa(\mu,\mu)$ is given on coadjoint orbits by 
	\begin{align}
		\mathbb P_\infty(\mu) = Z^{-1} e^{-\frac{2\theta}{\sigma^2} h(\mu)},  
	\label{MaxwellianDist}	
	\end{align}
	where $Z$ is the normalisation constant that enforces $\int \mathbb P_\infty(\mu) d\mu= 1$. 
\end{theorem}
Measures of the form \eqref{MaxwellianDist} are called \textit{Maxwellians}, or \textit{Gibbs measures}, for canonical ensembles in statistical physics. 
In statistical physics, the constant damping-to-forcing ratio ${2\theta}/{\sigma^2}$ would be associated with the inverse temperature $\beta= 1/(k_BT)$, where $k_B$ is the Boltzmann constant and $T$ is the Kelvin temperature. In this context, the normalisation constant $Z(\beta)$ is called the partition function, \cite{chirikjian2012stochastic2}.

\subsection{Random attractors}\label{RA-section}
The presence of both noise and dissipation in a dynamical system presents an interesting opportunity to study random attractors. 
We briefly describe the main steps in understanding these attracting sets, which we will follow explicitly in the example section, using numerical simulations.  

We refer the interested reader to \cite{crauel1994attractors,crauel1997random,arnold1995random,bonatti2006dynamics,kloeden2011nonautonomous} for extensive accounts of the topic of random attractors in the random dynamical systems theory.
Briefly, the invariant distribution $\mathbb P_\infty(\mu)$ of the Fokker-Planck equation represents the average solution of our dynamical system over all possible realisations of the noise, asymptotically in time, as $t\to \infty$. 
That is, averaging and taking the limit in time gives the probability measure, $\mathbb P_\infty$, which is usually smooth, covers the entire phase space and is independent of time and initial conditions. 
One can also take an alternative approach: instead of averaging over the ensemble of realisations of the noise, one may average only over the initial conditions, and let the system evolve toward large times. 
The distribution resulting from this procedure depends on time, and does not smoothly cover the entire phase space. 
Nevertheless, it can be shown that this distribution, called the random attractor, admits a probability measure, called the Sinai-Ruelle-Bowen measure (SRB). 
The SRB measure can be regarded as the equivalent of a volume preserving measure for non-volume preserving systems, e.g., when there is dissipation.  
The SRB measure plays an important role in the study and characterisation of chaotic dissipative dynamical systems.
Furthermore, the derivation of a general set of sufficient conditions for a dynamical system to admit a SRB measure is still an open problem. 
We refer to to the review paper \cite{young2002srb} for a detailed introduction to SRB measures and related open problems.
We will denote the SRB measure as $\mathbb{P}_\omega(\mu)$ for a given realisation of the noise $\omega$.
Under certain conditions, a strong link exists between $\mathbb P_\infty$ and the SRB measure $\mathbb P_\omega$, given formally by  
\begin{align}
	\int_\Omega \mathbb{P}_\omega(\mu) d\omega = \mathbb{P}_{\infty}(\mu),  
\end{align}
for the probability space $\Omega$; see \cite{crauel1998additive} for the derivation of these conditions.
The proof of this result is based on the observation that under double bracket dissipation the energy is monotonically decaying in time, and thus provides a Lyapunov function, which in turn implies the existence of attractive random sets. 
The stochastic process \eqref{SEP-Diss} does indeed admit random attractors, which may be singular sets. See \cite{schenk1998random,kondrashov2015data} and references therein for more details about this type of approach.

In this situation, following for example \cite{chekroun2011stochastic}, and provided that the largest Lyapunov exponent of the random system is positive (meaning that the system exhibits chaos), the existence of the non-singular SRB measure can be derived. 
Thus, the key step in establishing existence of an SRB measure is to determine a condition for the positivity of the largest Lyapunov exponent, as a function of the system parameters, especially the noise and dissipation amplitudes $(\sigma^2, \theta)$.
As in many examples of random attractors, the mathematical proof of positive top Lyapunov exponent is a nontrivial problem which involves evaluations of complicated integrals. 
Here, however, we will determine quantitative sufficient conditions for the positivity of the Lyapunov exponent via direct numerical simulations of the stochastic process and its linearisation.  
We refer the interested reader to \cite{arnaudon2016noise} for the numerical analysis of the stochastic rigid body or to \cite{lin2008shear, engel2016syncronisation} and references therein for other systems.

\section{Application with the stochastic free rigid body}\label{RB}
This section treats the classic example of the Euler-Poincar\'e dynamical equation; namely, the equation for free rigid body motion with three dimensions, described by the Lie group $SO(3)$.   
For a complete treatment from the viewpoint of reduction we refer to \cite{marsden1999intro}.
Of course, noise in the rigid body has already been considered in a number of previous works. (See for example  \cite{chirikjian2012stochastic1,chirikjian2012stochastic2} and references therein.) However, the system which we will obtain from this theory is quite different from those previously studied, as it preserves the geometry of the rigid body motion; in particular, it preserves the norm of the angular momentum. 

\subsection{The stochastic rigid body}
Before applying the theory outlined above, we should mention that we will use the isomorphism $\mathfrak{ so}(3)\cong \mathbb R^3$ which translates the commutator in the Lie algebra to the cross product of three-dimensional vectors, via $[A , B ] \to \boldsymbol A\times\boldsymbol B$, where vectors in $\mathbb R^3$  are denoted with bold font.
This allows us to use the scalar product of vectors as our pairing, via the formula $\boldsymbol A\cdot \boldsymbol B= -\frac12 \kappa(A,B)$. 
We skip the details and directly use the reduced Lagrangian of the free rigid body 
\begin{align}
	l(\boldsymbol \Omega) = \frac12 \boldsymbol \Omega \cdot\mathbb{ I} \boldsymbol \Omega
	: = \frac12 \boldsymbol \Omega\cdot \boldsymbol \Pi\,,
	\label{H-RB}
\end{align}
where $\boldsymbol \Omega$ is the angular velocity, $\mathbb I= \mathrm{diag}(I_1,I_2,I_3)$ is a prescribed moment of inertia and  $\boldsymbol \Pi$ is the angular momentum.
Notice that the Legendre transform gives the reduced Hamiltonian $h(\boldsymbol \Pi) = \frac12 \boldsymbol \Pi \cdot \,\mathbb I ^{-1}\boldsymbol \Pi$.
We take the stochastic potential to be linear in the momentum variable $\boldsymbol \Pi$
\begin{align}
	\Phi_i(\boldsymbol \Pi) = \sum_{i=0}^3 \boldsymbol{\sigma}_i\cdot \boldsymbol\Pi\,,
	\label{Phi-RB}
\end{align}
where the constant vectors $\boldsymbol{\sigma}_i$ span $\mathbb R^3$.
The stochastic process for $\boldsymbol \Pi$ is then computed from \eqref{SEP} to be 
\begin{align}
	d\boldsymbol\Pi + \boldsymbol\Pi\times \boldsymbol\Omega\, dt + \sum_i\boldsymbol\Pi\times \boldsymbol{\sigma}_i \circ dW^i_t=0.
	\label{Sto-RB-stra}
\end{align}
One can check that for either the It\^o or Stratonovich stochastic process the Casimir level set $\|\boldsymbol\Pi\|^2=\mathrm{c}^2$ is preserved even with the noise. The Casimir level set defines the momentum sphere of radius $c$, which is the coadjoint orbit,
Although the energy $h(\boldsymbol\Pi)$ is not a conserved quantity, it stays bounded within the maximum and minimum energies of the deterministic system, as the dynamics takes place on the momentum sphere, \cite{arnaudon2016noise}.  
The invariant solution of the associated Fokker-Planck equation (see below in \eqref{FP-RB-SD} with $\theta=0$) is a constant on the coadjoint orbit, or momentum sphere. 
Thus, this system behaviour is similar to that of the heat equation on a compact domain, but with the non-trivial geometry of the Casimir level set. 

\subsection{Adding dissipation}
The double bracket dissipation for the rigid body involves the only Casimir $\|\boldsymbol \Pi\|^2$ and gives the dissipative stochastic process in \eqref{SEP-Diss}, 
\begin{align}
	d\boldsymbol\Pi + \boldsymbol\Pi\times \boldsymbol\Omega\,  dt + \theta\, \boldsymbol \Pi\times( \boldsymbol \Pi\times \boldsymbol \Omega)\, dt+ \sum_i\boldsymbol\Pi\times \boldsymbol{\sigma}_i \circ dW^i_t=0 .
	\label{RB-diss}
\end{align}
In the absence of noise, the energy decay of the deterministic dissipative rigid body is given by,
\begin{align}
	\frac{d h}{dt}= -\theta \|\boldsymbol \Pi\times \boldsymbol \Omega \|^2\,.
	\label{energy-decay}
\end{align}
Consequently, the dissipation will bring the system to one of the minimal energy positions, where $\Pi$ and $\Omega$ are aligned, which corresponds to (relative) equilibria. 

In the presence of the noise, the associated the Fokker-Planck equation may be found as
\begin{align}
	\begin{split}
	\frac{d}{dt}\mathbb P& + (\boldsymbol\Pi\times \boldsymbol\Omega)\cdot  \left ( \nabla \mathbb P - \theta\, \boldsymbol \Pi \times  \nabla \mathbb P\right )+\frac12\sum_i (\boldsymbol\Pi\times \boldsymbol{\sigma}_i)\cdot \nabla[(\boldsymbol\Pi\times \boldsymbol \sigma_i) \cdot \nabla \mathbb P ]=0.
	\end{split}
	\label{FP-RB-SD}
\end{align}
Unlike the case with $\theta=0$, this equation will not have a constant invariant solution. 
As an illustration, we derive the invariant distribution $\mathbb P_\infty$ of Theorem \ref{FP-diss-thm} for this simple case. 
First, we rewrite the Fokker-Planck equation \eqref{FP-RB-SD} as 
\begin{align}
	\frac{d}{dt}\mathbb P + (\boldsymbol\Pi\times \boldsymbol\Omega) \cdot \nabla \mathbb P + \nabla\cdot \left (  \theta\, \boldsymbol \Pi \times (\boldsymbol\Pi\times \boldsymbol\Omega) \mathbb P   -\frac12\sigma^2\, \boldsymbol\Pi\times( \boldsymbol\Pi\times \nabla \mathbb P )\right )=0 \,,
\end{align}
where we have used $\nabla\cdot (\boldsymbol\Pi\times (\boldsymbol\Pi\times \boldsymbol\Omega))= 0$. The last term in \eqref{FP-RB-SD} simplifies as
\begin{align*}
	\sum_i (\boldsymbol\Pi\times \boldsymbol \sigma_i) [( \boldsymbol\Pi\times  \boldsymbol \sigma_i)\cdot \nabla\mathbb P] 
	&=  \sum_i (\boldsymbol\Pi\times \boldsymbol \sigma_i) [( \nabla\mathbb P \times \boldsymbol\Pi)\cdot  \boldsymbol \sigma_i] 
	= \boldsymbol\Pi\times (\nabla\mathbb P\times \boldsymbol\Pi), 
\end{align*}
since the sum over $i$ is simply the decomposition of the vector $( \nabla\mathbb P \times \boldsymbol\Pi)$ into its $\boldsymbol \sigma_i$ basis components. 
Hence, the invariant solution is the Gibbs measure, given by
\begin{align}
	\mathbb P_\infty(\boldsymbol \Pi) = Z^{-1}e^{-\frac{2\theta}{\sigma^2}h(\boldsymbol \Pi)}\,.
\end{align}
We recover the constant solution when $\theta=0$. Notice that when $\sigma=0$, the equilibrium distribution $\mathbb P_\infty$ has singular support, which comprises two Dirac delta functions at the lowest energy relative equilibrium points, as expected from the double bracket dissipation.   

\subsection{Random attractors}

We will end this short note by commenting on the random attractor of the stochastic rigid body. 
The proof of its existence was sketched earlier in the paper, and can be found in more detail in \cite{arnaudon2016noise}. 

The existence of random attractive sets is a direct consequence of the dissipation, although showing that this set is non-singular requires a certain amount of care. 
Indeed, chaotic motion only occurs, provided the noise amplitude is sufficiently large, compared to the dissipation. 
Chaotic motions are characterised as possessing positive Lyapunov exponents, which describe the sensitivity of a dynamical system to the initial conditions. In particular, chaos via sensitive dependence on initial conditions occurs when the largest Lyapunov exponent is positive. This condition may be achieved for the stochastic rigid body, by choosing appropriate dissipation and noise coefficients. 
We refer to \cite{arnaudon2016noise} for a deeper analysis of the top Lyapunov exponent, and simply highlight its key property here.
Namely, positivity of the top Lyapunov exponent requires sufficient shear as well as strong enough noise.
In the rigid body, the notion of shear is given by the fact that nearby orbits have different speeds, characterised by the relative values of the moment of inertia, and the absolute value of the momentum.  
This shear is a key ingredient, as it affords the formation of the stretching and folding mechanism described below. 
This type of random attractor generated by shear is well known for motions in the plane, and has been studied, for example, in \cite{lin2008shear,engel2016syncronisation}.

\begin{figure}[h]
	\centering
	\subfigure{\includegraphics[scale=0.4]{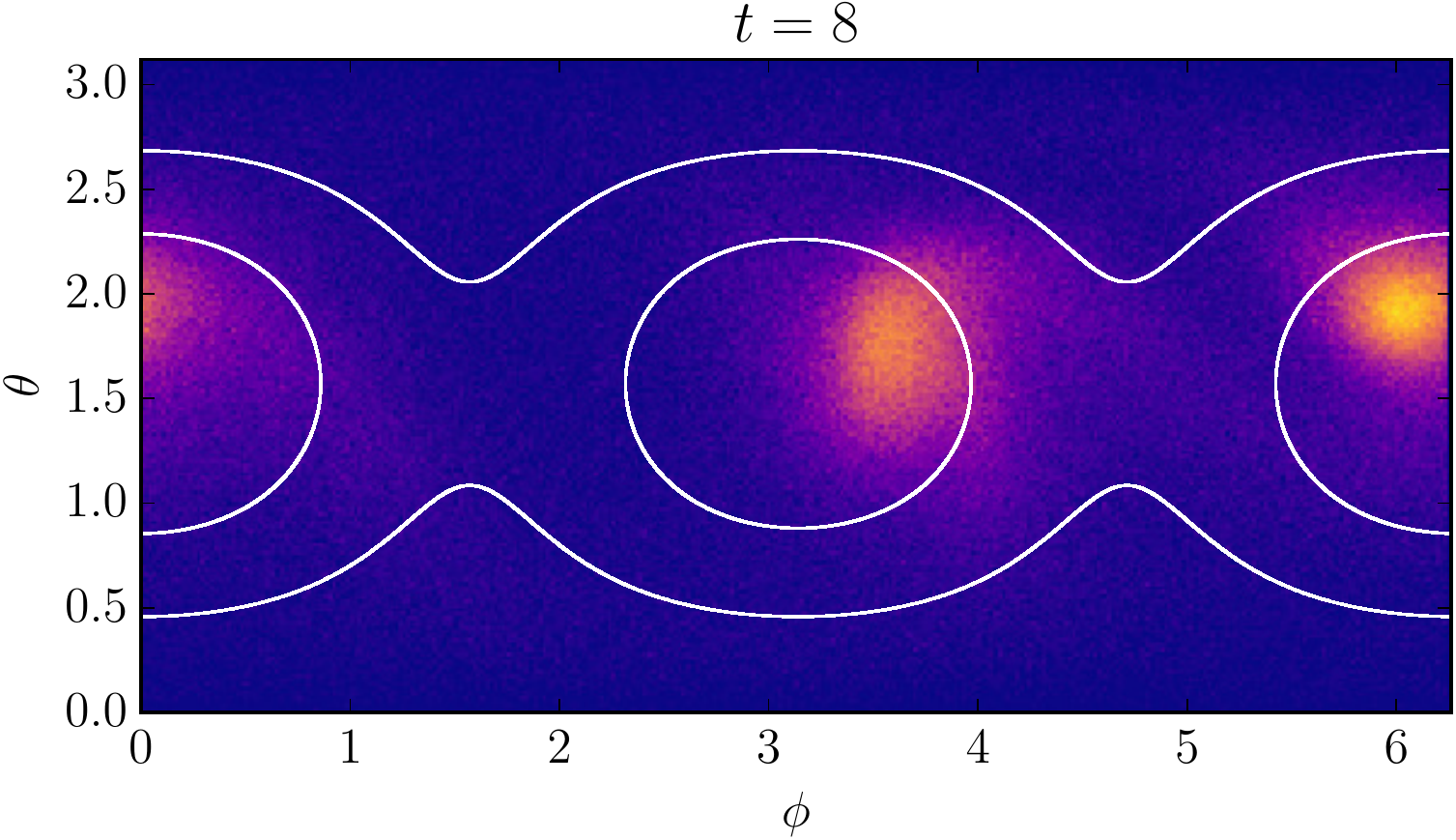}}
	\subfigure{\includegraphics[scale=0.4]{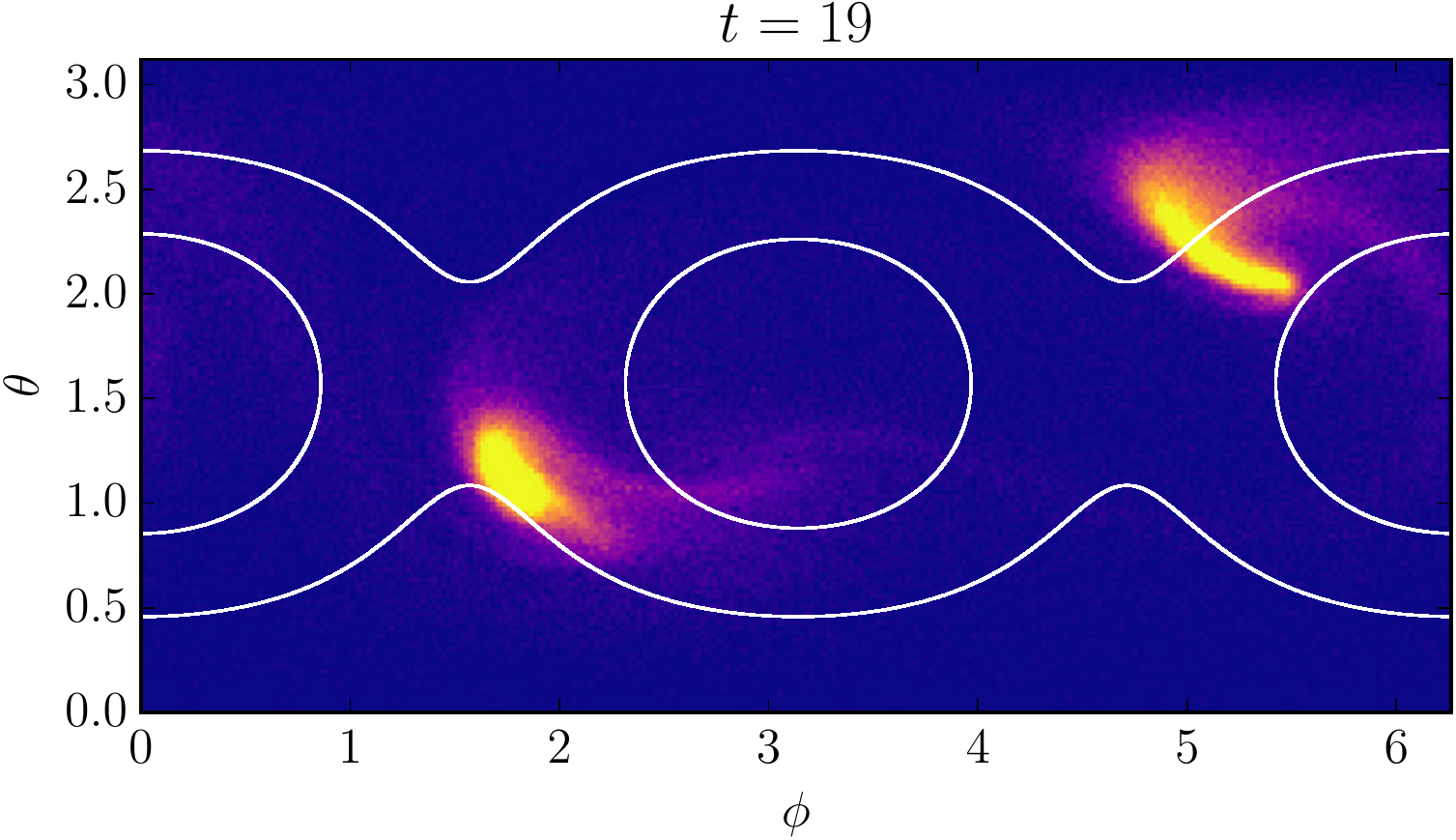}}
	\subfigure{\includegraphics[scale=0.4]{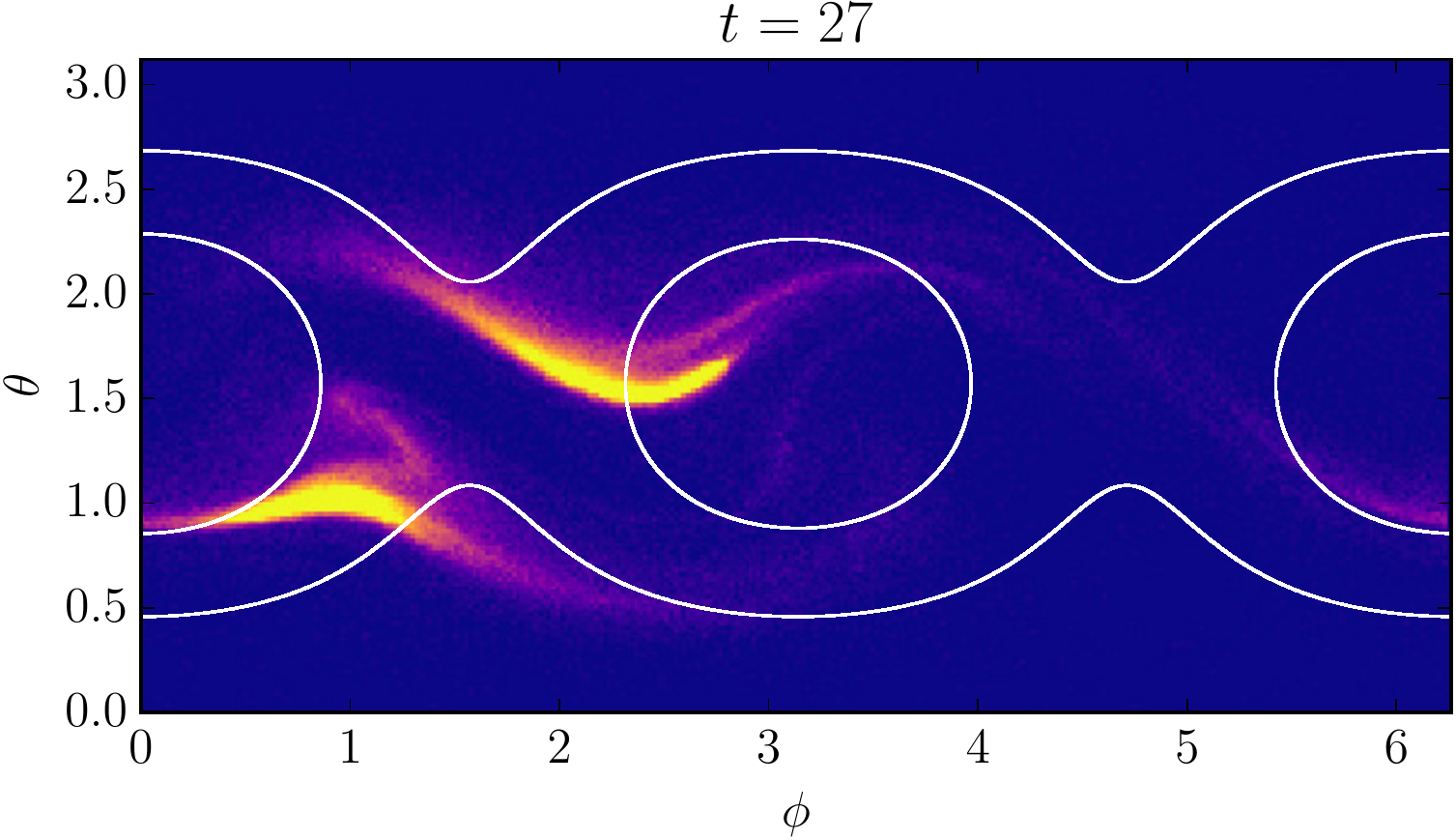}}
	\subfigure{\includegraphics[scale=0.4]{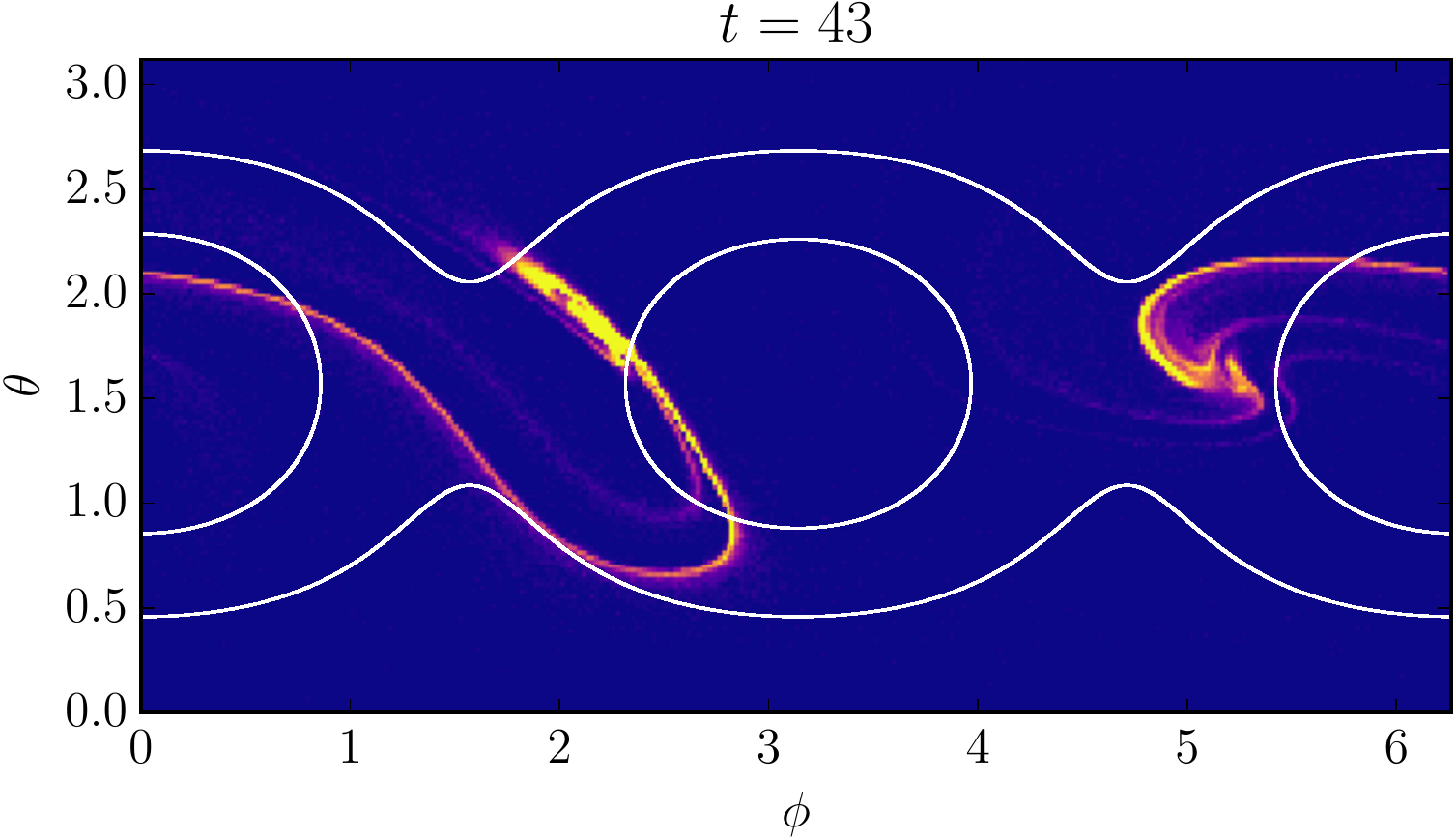}}
	\caption{{\small We display four snapshots of the same rigid body random attractor with $\mathbb I= \mathrm{diag}(1,2,3)$, $\theta=0.5$ and $\sigma= 0.5$. 
	The simulation started from a uniform distribution of rigid bodies on the momentum sphere at $t=0$. The color is in log scale and we simulated 400,000 rigid body initial conditions with a split step numerical scheme. }}
	\label{fig:RB-RA}
\end{figure}

Having shown that non-singular random attractors exist on coadjoint orbits for the stochastic rigid body with selective decay, one may hope to characterise the properties of these attractors.
At present, we have made little progress in this interesting endeavour. However, some clues to these properties may be found in numerical simulations. 
For example, we display in Figure \ref{fig:RB-RA} a realisation of a random attractor of the rigid body.
What is shown is the probability density, in log scale, calculated from a Monte-Carlo simulation of the stochastic rigid body equations. 
This probability density approximates the SRB measure, which is supported on the random attractor, for a certain noise and dissipation. 
This measure is time dependent, and its motion exhibits stretching and folding, which is a common feature of attractors with both positive and negative Lyapunov exponents. 
Indeed, the positive exponent produces stretching of this set, and the negative one produces compression, which together with the original nonlinear deterministic rigid body dynamics, creates the folding process.  
In principle, repeated iterations of this stretching and folding mechanism could create a structure similar to a horseshoe map, and this map would prove sensitivity to initial conditions, although we have not been able to prove that this construction exists in the presence of stochasticity.  

\subsection*{Acknowledgements}
{\footnotesize
We are grateful to many people for fruitful and encouraging discussions, including S. Albeverio, J.-M. Bismut, N. Bou-Rabee, M. D. Chekroun, G. Chirikjian, D. O. Crisan, A.-B. Cruzeiro, J. Eldering, M. Engel, N. Grandchamps,  P. Lynch, J.-P. Ortega, G. Pavliotis, V. Putkaradze, T. Ratiu and C. Tronci.
The simulations were run with the Imperial College High Performance Computing Service.
We also acknowledge the Bernoulli Centre for Advanced Studies at EPFL where parts of this work were elaborated.  
AA acknowledges partial support from an Imperial College London Roth Award and AC from a CAPES Research Award BEX 11784-13-0. 
All the authors are also supported by the European Research Council Advanced Grant 267382 FCCA held by DH.
Last, but not least, we want to mention the inspiration of the lectures of H. Dumpty on broken symmetry, and stochastic processes on coset spaces. 
}

\bibliographystyle{alpha}
\bibliography{sgm-2015}

\end{document}